\documentclass{amsart}

\usepackage{amsthm,amssymb,verbatim}
\usepackage{graphicx}
\usepackage{enumerate}

 \newcommand{\MM}{\mathcal{M}}

 \newcommand{\FF}{\mathcal{F}}
 \newcommand{\RR}{\mathbf{R}}  
 \newcommand{\ZZ}{\mathbf{Z}}  
 \newcommand{\BB}{\mathbf{B}}  
 \newcommand{\dist}{\operatorname{dist}}

 \newcommand{\eps}{\epsilon}

\input epsf
\def\begfig {
\begin{figure}
\small }
\def\endfig {
\normalsize
\end{figure}
}

    \newtheorem{claim}{Claim} 
    \newtheorem*{claim*}{Claim}
    \newtheorem*{theorem*}{Theorem}

    \newtheorem{theorem}    {Theorem}     

    \theoremstyle{definition}

    \theoremstyle{definition}
    
    \newtheorem*{remark*}{Remark}

\usepackage{hyperref}
\usepackage[alphabetic, msc-links]{amsrefs}

\begin{document}

\setlength{\baselineskip}{1.2\baselineskip}

\renewcommand{\thesubsection}{\thetheorem}

\title[Sequences of Embedded Minimal Disks]{Sequences of embedded minimal disks
   whose curvatures blow up on a prescribed subset of a line}
\author{David Hoffman}
\address{Department of Mathematics\\ Stanford University\\ Stanford, CA 94305}
\email{hoffman@math.stanford.edu}
\author{Brian White}
\thanks{The research of the second author was supported by the NSF
  under grant~DMS-0104049}
\email{white@math.stanford.edu}
\date{November 13, 2009; revised August 31, 2011}


\begin{abstract}
For any prescribed closed subset of a line segment in Euclidean 3-space, we construct a sequence of minimal disks that are properly embedded in an open solid cylinder around that line and  that have
curvatures blowing up precisely at the points of the closed set. 
\end{abstract}

\maketitle

\section*{Introduction}\label{section:intro}
\newcommand{\LL}{\mathcal{L}}     
\newcommand{\HH}{\operatorname{\mathcal H}}
\newcommand{\HJ}{\operatorname{\mathcal H}_J}
\newcommand{\DD}{\operatorname{\bf\mathcal D}}   
\newcommand{\CylJ}{\BB(0,1)\times z(J)}                    

 In this paper, we prove the following result about convergence
 of properly embedded minimal disks to a lamination:

 \begin{theorem}\label{maintheorem}
 
 Let
 \[
       C=\BB(0,1)\times (a,b)\subset \RR^3
\]
  be an open solid circular cylinder, possibly
 infinite, that is rotationally symmetric about the $z$-axis, $Z$.
 Let $K$ be a relatively closed subset of $Z\cap C$.
Then there exists a sequence of minimal disks that are properly embedded in $C$ and that have the
following properties:
\begin{enumerate}[\upshape (1)]
 \item\label{Statement-1}
The curvatures of the disks are uniformly bounded on
compact subsets of $C \setminus K$, and the curvatures blow up  at each point of $K$.\footnote{The
curvatures of a sequence of  minimal disks $D_n$ blow up at a point
$p$ if there exists a sequence $\{p_n\in D_n\}$ with $p_n\rightarrow
p$ such that the absolute value of the
curvature of $D_n$ at $p_n$ tends to infinity as $n\rightarrow \infty$}

  \item\label{Statement-2} The  minimal disks converge to a limit lamination of $C\setminus K$ consisting
  of the following leaves:
 \begin{enumerate}[\upshape (i)]
\item\label{Statement:flat-leaves} For each $p\in K$, the horizontal punctured unit disk centered at $p$:
\[
   \DD_p = (C\setminus Z) \cap \{ z = z(p)\};
\]
   
\item\label{Statement:other-leaves}  For each component $J$ of  $(Z\setminus K)\cap C$, a leaf $\MM_J $ that is 
 properly embedded in the cylinder
 \[
       \CylJ.
  \]
 The leaf $\MM_J$ contains the segment $J$ and is therefore symmetric under rotation
 by $180^{o}$ around $Z$.  Each of the two connected components 
 of $\MM_J\setminus J$ is an infinitely-sheeted multigraph over the punctured unit disk, 
 and $\partial z/\partial \theta>0$ everywhere on $\MM_J\setminus J$.
\end{enumerate}
\item\label{Statement-3}  The lamination 
 extends smoothly to a lamination of $C\setminus \partial K$,
 but it does not extend smoothly to any point in $C\cap \partial K$.
 \end{enumerate}
 \end{theorem}
 
 Here $z(J)$ denotes the image of $J$ under the coordinate map $z: \RR^3\to \RR$. 
 
 Note that according to Statement~\ref{Statement:other-leaves}, each of the components of $\MM_J\setminus J$
 can be parametrized as
 \[
    (r,\theta)\in (0,1) \times \RR \mapsto (r\cos\theta, r\sin\theta, f(r,\theta)).
\]
Furthermore, if we let $(c,d)=z(J)$, then
(by properness)
\begin{align*}
     \lim_{\theta\to-\infty} f(r,\theta) &= c,  \\
     \lim_{\theta\to \infty} f(r,\theta) &= d,
 \end{align*}
 the convergence being uniform away from $r=0$.

Theorem~\ref{maintheorem} is well-known in  case $K$ is the entire interval $Z\cap C$:  let the $n$th disk
be the standard helicoid scaled by a factor of $1/n$ and restricted to $C$.  
Theorem~\ref{maintheorem} was proved 
by Colding and Minicozzi \cite{ColdingMinicozziExample} when $K$ consists of  a single point,
by Brian Dean~\cite{BrianDean} when $K$ is an arbitrary finite set of points, and 
very recently by Siddique Kahn~\cite{SiddiqueKahn}
when $K$ is an interval with exactly one endpoint in $C$.  
Our work was inspired by Kahn's result, although the methods are very different:
Kahn, like Colding, Minicozzi, and Dean, used the Weierstrass Representation, whereas our
approach is variational.  Stephen~J.~Kleene has given a different proof of
Theorem~\ref{maintheorem} using the Weierstrass representation; see~\cite{Kleene}.

In the study of minimal varieties, the known examples of singularities have been rather tame.  In particular, we believe
that Theorem~\ref{maintheorem} provides the first examples of Cantor sets of singularities and of singular sets with
non-integer Hausdorff dimension.

Our  paper is organized as follows. 
In  section~\ref{context}, we discuss what is known in general
about curvature blow up in sequences of properly
embedded minimal disks.
We prove Theorem~\ref{maintheorem}
in section~\ref{mainsection}. The proof depends on results in
Sections~\ref{ExistenceUniqueness} and ~\ref{RadoCurvatureSection}.
In Section~\ref{ExistenceUniqueness}, we prove existence and
uniqueness theorems for embedded minimal disks with certain
rotationally symmetric boundaries.  In
Section~\ref{RadoCurvatureSection}, we use Rado's Theorem to deduce
curvature estimates for our examples.

\section{The General Context}\label{context}

Our paper gives examples of curvature blow up in sequences of properly embedded
minimal disks.  In this section, we describe what is known in general about
such curvature blow up.  (The results described in this section are
not used in the rest of the paper.)

Suppose that for each $n$ we have a minimal disk $D_n$ that is  properly embedded in an
open subset $U_n$ of $\RR^3$, where $U_n  \subset U_{n+1}$ for each $n$.  Let $U=\cup_nU_n$.
By passing to a subsequence, we may suppose that there is a relatively closed subset $K$ of $U$
such that the curvatures of the $D_n$ blow up at all points of $K$ and such that the $D_n$ converge
smoothly in $U\setminus K$ to a limit lamination $\LL$ of $U\setminus K$.
It is natural to ask how general the set $K$ can be, and to what extent the lamination and/or the leaves
of the lamination can be smoothly extended to include points in~$K$.
In particular,
\begin{enumerate}[Q1.]
\item If $p$ is a point in the blow-up set $K$, 
must there be a leaf of the lamination that extends smoothly across~$p$?
\item Must $K$ locally be contained in a nice (e.g., $C^1$ or perhaps $C^{1,1}$) curve? 
\item Given an arbitrary closed subset $S$ of a nice 
 curve 
 in the open set $U$, is there an example for which which the blow-up set $K$
 is precisely~$S$?
\end{enumerate}

(More precisely, Q2 is: must each point $p\in K$ have a neighorhood $W$ such that $K\cap W$
is, after a rotation, contained in the graph of a $C^1$ (or $C^{1,1}$) function
from $\RR$ to $\RR^2$?)

According to Theorem~5.8 of~\cite{ColdingMinicozziII}, the answer to question~Q1
is ``yes''.  
In particular, every point of $K$ has a neighborhood $W$ such that if $p\in  K\cap W$, then
there is a leaf $L_p$ of the limit lamination such that $L'_p:=(L_p\cup\{p\})\cap W$ is a properly embedded
minimal submanifold of $W$.  Furthermore, if we choose $W$ small enough, then the $L'_p$ are
all graphs over a common plane (with uniformly small slopes.)  Note that for each $p\in K\cap W$, there
is only one such $L'_p$ by the strong maximum principle.

The answer to Q2 is also ``yes", at least in the $C^1$ case.  That the answer is ``yes"
for Lipschitz curves is implicit in the work of Colding and Minicozzi;  
see, for example,~\cite{ColdingMinicozziIV}*{Section I.1} and~\cite{ColdingMinicozziIV}*{Theorem 0.1} for a very
similar result.  
Assuming the results of Colding and Minicozzi and an improvement by Meeks~\cite{MeeksRegularity},
\cite{WhiteC1} proves that the answer to Q2 is ``yes'' in the $C^1$ case.

Colding and Minicozzi also proved
 that if $U=\RR^3$ (the so-called ``global case'') and if $K$ is nonempty,
  then (after a rotation) $K$ is the graph of a lipschitz
function $z\in \RR \mapsto (x(z), y(z))$
and the lamination is the foliation by horizontal planes punctured 
at the points of $K$ \cite{ColdingMinicozziIV}*{Theorem 0.1}.  
 In particular, the lamination
extends smoothly to all of $\RR^3$.   
By~\cite{MeeksRegularity} (described below), the curve is in fact a straight line that is perpendicular
to the planes.

In the local case $U\ne \RR^3$, the behavior can be very different, as  Theorem~\ref{maintheorem} indicates.
Meeks proved that if the lamination extends smoothly to a foliation of $U$ and if $K$ is a Lipschitz curve
that intersects the leaves transversely, then $K$ is a $C^{1,1}$ curve 
and it intersects the leaves orthogonally~\cite{MeeksRegularity}.
Meeks and Weber~\cite{MeeksWeber} constructed an example
for which $U= \RR^3\setminus Z$, the blow-up set $K$ is a horizontal circle centered at a point in $Z$,
and the limit lamination consists of the vertical half-planes with $Z$ as edge, punctured by $K$.
They go on to prove that, given any $C^{1,1}$ curve, there is an example in which the blow-up set  $K$ is that curve, 
$U$ is tubular neighborhood of $K$, and the lamination is a foliation of $U$ by planar punctured disks orthogonal
to $K$.  The examples of Meeks and Weber and the examples in this paper suggest that
the answer to Q3 may be ``yes".

\section{The proof of  Theorem ~\ref{maintheorem}}\label{mainsection}

\newcommand{\cyl}{\operatorname{Cyl}}
\newcommand{\interior}{\operatorname{Interior}}

\begin{proof}

It suffices to prove the theorem for the unbounded cylinder $C=\BB(0,1)\times \RR$:
the case of bounded cylinders follows by restriction.

For each connected component $J$ of  $Z \setminus K$, choose a
smooth embedded curve $S_J$ in $\partial \BB(0,1)\times J$ such that (see figure 1)
\begin{enumerate}[ \upshape (1)]
\item The projection $(x,y,z)\mapsto (0,0,z)$ induces a diffeomorphism from $S_J$ to $J$.
\item The derivative $d\theta/dz$ is strictly positive at each point of $S_J$. (Here $\theta$ is the
angle of the cylindrical  coordinate of a point in $\RR^3\setminus Z$.
Of course, $\theta$ is defined only up to an integer multiple of
$2\pi$, but $d\theta/dz$ is well defined.)
\item The curve $S_J$  winds around the cylinder
infinitely many times as $z\rightarrow c$ and as $z\rightarrow d$, where $c$ and $d$ are the infimum and
supremum, respectively, of $z$ on~$J$.   
In other words, $\lim_{z\to c}\theta(z)= -\infty$ and $\lim_{z\to d}\theta(z)=\infty$.
\end{enumerate}

 \begfig
 \hspace{0.3in}
  \vspace{.2in} \centerline {
  \includegraphics[width=1.25in]{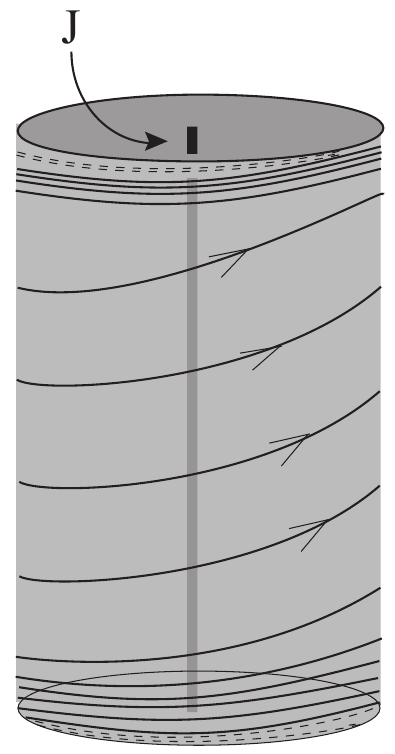}
   }
 \hfill
 \begin{center}
 \vspace{0.3cm}
 \parbox{4.1in}{
  \caption{\label{figure:genus-one helicoids} 
The curve $S_J$ on  $\partial \BB (0,1)\times \RR$, as described in conditions (1), (2), and (3). Here, $J$ is a bounded open interval of the axis $Z$, the endpoints of which are in the closed set $K$. The limit set of  $S_J$ consists of the two horizontal circles on $\partial \BB (0,1)\times \RR$ at the heights of the endpoints of $J$.  }
 }
 \end{center}
 \endfig

Now choose numbers $a_n\to -\infty$ and $b_n\to\infty$, and choose
choose smooth curves $\gamma_n$ in $\partial \BB(0,1)\times [a_n,b_n]$
such that
\begin{enumerate}[(i)]
\item The curve $\gamma_n$ can be smoothly parametrized by $z\in [a_n,b_n]$,
(That is, $(x,y,z)\mapsto z$ induces a diffeomorphism
  from $\gamma_n$ to $[a_n,b_n]$.)

\item The derivative $d\theta/dz$ is strictly positive at each point of $\gamma_n$.

 \item The $\gamma_n$ converge
smoothly to the lamination of $\partial C$ consisting of the $S_J$'s together with
the horizontal circles of radius $1$ centered at points of $K$.
\newcounter{enumi-saved}
\setcounter{enumi-saved}{\value{enumi}}
\end{enumerate}

Let $\Gamma_n$  be the simple closed curve consisting of $\gamma_n$
and $Z\cap \{a_n \le z\le b_n\}$, together with two radial segments at heights 
in $z=a_n$ and $z=b_n$. By
Theorem~\ref{existence2} in section~\ref{ExistenceUniqueness}, $\Gamma_n$ bounds a unique  embedded minimal
disk $D_n$. 

Extend $D_n$ by repeated Schwartz reflection in the the top and bottom edges to get
an infinite minimal strip $\tilde D_n$.   The boundary of $\tilde D_n$ has two
components: the axis $Z$, and the helix-like curve $\tilde \gamma_n$ in $\partial C$ obtained from $\gamma_n$ by 
iterated reflection. 

To ensure that the boundary curve $\tilde \gamma_n$ is smooth, we
need to impose one additional condition on the choice of $\gamma_n$:
\begin{enumerate}[(i)]
\setcounter{enumi}{\value{enumi-saved}}
\item\label{even-order} The even-order derivatives $d^k\theta/dz^k$ ($k=2,4,\dots$) vanish at the endpoints of $\gamma_n$.
\end{enumerate}
(If this condition were not satisfied, then the curve $\tilde \gamma_n$ would not be smooth
at the endpoints of $\gamma_n$.)

By Theorem~\ref{existence2}, $\tilde D_n$ and its rotated images $R_\theta \tilde D_n$ foliate
 $C\setminus Z$.    Thus if $D_n^* = \tilde D_n \cup \rho_Z \tilde D_n \cup Z$
is the minimal surface obtained from $\tilde D_n$ by Schwartz reflection in $Z$, then $D_n^*$
is embedded.  Here $\rho_Z=R_\pi$ denotes rotation by $\pi$ about $Z$.

The disks $D_n^*$ are the disks whose existence is asserted in the statement of the theorem.
We now prove that they have the properties asserted by the theorem.

By Theorem~\ref{existence2} applied to $D_n$, each of the two connected components of
$D_n^*\setminus Z$  is a multigraph over $\BB(0,1)\setminus \{0\}$.  (Of course those components are
 $\tilde D_n$ and $\rho_Z\tilde D_n$.)
Thus if $H$ is an open halfspace bounded by a plane containing $Z$, then
$H\cap D_n^*$ is a union of graphs over a half-disk. Standard
estimates for solutions of the minimal surface equation then imply
that the principal curvatures of the closures of the $D_n^*$ are uniformly bounded on
compact subsets of $H$. Thus after passing to a subsequence, those
graphs will converge smoothly (away from $\partial H$) to a collection $\mathcal{G}_H$ of
minimal graphs over a half disk.   Note that if $G$ is such a graph, then
$(\partial G)\cap H$ must be one of the following:
\begin{enumerate}[\upshape (H1)]
\item\label{H1} a horizontal semicircle at height $z(p)$ for some $p\in K$, or
\item\label{H2} one of the connected components of $S_J\cap H$ or one of the connected
components of $(\rho_Z S_J)\cap H$, where $J$
 is a connected component of $Z\setminus K$.  
 \newcounter{enumi-H-saved}
 \setcounter{enumi-H-saved}{\value{enumi}}
\end{enumerate}
Furthermore,
\begin{enumerate}[\upshape (H1)]
\setcounter{enumi}{\value{enumi-H-saved}}
\item\label{H3} If $\Gamma$ is a connected component of $S_J\cap H$ or of $(\rho_Z S_J)\cap H$, 
 then there is exactly one
 graph $G$ in $\mathcal{G}_H$ whose boundary (in $H$) is $\Gamma$.
\end{enumerate}

Since $H$ is arbitrary, this means that (after passing to a
subsequence) the $D_n^*\setminus Z$ will converge smoothly on
compact subsets of $\RR^3\setminus Z$ to a lamination
$\LL$ of $C\setminus Z$ consisting of a union of multigraphs.  
 Because each leaf is embedded, it must either have 
 a single sheet (and thus be a graph) or else have infinitely many sheets.

Let $J$ be a component of $Z\setminus K$.  
  By Theorem~\ref{RadoCurvEst},
principal curvatures of the $D_n^*$ are uniformly 
bounded\footnote{We could apply standard curvature estimates
                instead of Theorem~\ref{RadoCurvEst} if we had uniform area 
                bounds in little balls centered along $Z$, but proving the necessary
                area bounds does not seem to be any easier than proving 
                Theorem~\ref{RadoCurvEst}.} 
on compact subsets of $ \BB(0,1)\times z(J)$.  Since the principal
 curvatures are also uniformly bounded on compact subsets of  $\RR^3\setminus Z$,  this means that the
curvatures are in fact uniformly bounded on compact subsets of $\RR^3\setminus K$. 
Thus (perhaps after passing to a further subsequence) the $D_n^*$ converge
smoothly on compact subsets of $C\setminus K$ to a lamination $\LL'$ of that region. 
 Of course $\LL$ is the restriction of $\LL'$ to the gutted cylinder $C \setminus Z$.

\begin{claim}\label{FirstClaim}
  The horizontal circle ${\mathcal C}_p$ of radius $1$
centered at a point  $p\in K$ bounds a unique leaf of ${\mathcal
D}_p\in \LL'$.  That leaf is the planar punctured disk bounded by ${\mathcal C}_p$.
\end{claim}

\begin{proof}[Proof of Claim~\ref{FirstClaim}.]  As above, let $H$ be an open halfspace of
$\RR^3$ with $Z\subset \partial H$.   Then $\LL'\cap H$ is a
union of minimal graphs over a half-disk.   Consider
the set $Q$ of those graphs that contain ${\mathcal C}_p\cap H$ in
their boundaries. That set is compact, so there is an uppermost
graph $G_H$ in $Q$. Note that the union of the $G_H$ as $H$ varies
(by rotating it around $Z$) forms a single smooth minimal graph $G$
over the punctured disk $\BB(0,1)\setminus \{0\}$. That graph satisfies
the minimal surface equation. As a minimal surface in $\RR^3$, the
boundary of $G$ is the circle ${\mathcal C }_p$ together with some or
all of $Z$. Since a solution to the minimal surface equation cannot
have an isolated interior singularity, the graph extends to a regular
minimal surface over $\BB(0,1)$.  (See Theorem 10.2 of 
\cite{OssermanBook}. The result is originally due to Bers \cite{Bers}, but  the
proof in \cite{OssermanBook} using catenoidal barriers is due to Finn \cite{Finn}.)
 Since the boundary values define a planar circle, the
graph must be a flat disk, so $G$ is a flat planar disk. Recall that
$G$ is the uppermost leaf of $\LL$ that contains ${\mathcal C}_p$. By
the same argument, it is also the lowermost leaf of $\LL$ containing
${\mathcal C}_p$. Thus it is the unique leaf in $\LL$
containing ${\mathcal C}_p$. 
\end{proof}

\begin{claim}\label{2ndClaim}  If $p\in K$, then there is a sequence $p_n\in
D_n^*$ converging to $p$ such that  the norm of the second fundamental form of $D_n^*$ at $p_n$
tends to infinity.
\end{claim}

\begin{proof}[Proof of Claim ~\ref{2ndClaim}.] Suppose not.  Then there is a ball $\BB\subset \RR ^3$
centered at $p$ and a subsequence of the $D_n^*$ (which we may take
to be the original sequence) such that the curvatures of the $D_n^*$
are uniformly bounded in $\BB$. It follows that the lamination
$\LL'$ extends smoothly to a lamination
of $\BB$ and that  the convergence of
$D_n^*\cap \BB$ to the lamination of $\BB$ is smooth. By Claim 1,
 the leaf 
containing $p$ is a horizontal disk. But each
$D_n^*$ contains the axis $Z$, and therefore  has a vertical tangent
plane at $p$. Hence  the convergence cannot be smooth. This
contradiction proves the claim.
\end{proof}

We have now completely established Statement~\ref{Statement-1} of the Theorem~\ref{maintheorem},  
and we have established that there is a limit lamination of
 $C\setminus K$. We also know that  for each $p\in K$, the punctured disk $\DD_p$
 is a leaf of  this lamination. Thus we have proved 
 Statement~\ref{Statement:flat-leaves}.

To prove Statement~\ref{Statement:other-leaves}, let $J= \{ p\in Z: c<z(p)< d \}$ be one of the components of $Z\setminus K$.
We now analyze the leaves of the foliation that lie between the punctured disks at heights $c$ and $d$.
Let $\MM_J^+$ and $\MM_J^-$ be the leaves of $\LL$ that contain $S_J$ and $\rho_Z(S_J)$,
respectively.  Note that by~(H\ref{H1}) and~(H\ref{H2}),
these are both infinite covers of $\BB(0,1)\setminus \{0\}$.
By~(H\ref{H3}), there are no other leaves of $\LL$ in the region $\{c<z<d\}$.

Note that $J\subset D_n^*$ for every $n$, so $J$ is contained in one of the leaves $\MM_J$
of the lamination $\LL'$.   Now $\MM_J$ is simply connected since each $D_n^*$ is
simply connected,
so $\MM_J\setminus Z$ must contain two components. These components are leaves of $\LL$, so they
must be $\MM_J^+$ and $\MM_J^-$.   Thus $\MM_J = J\cup \MM_J^+ \cup \MM_J^-$ is the unique leaf of $\LL'$ in 
the slab $\{ c < z < d\}$.   

Let $H$ be an open halfspace of $\RR^3$ bounded by a plane containing $Z$.  
Note that the components of $\MM_J\cap H$ form a countable discrete set corresponding to the countable
discrete set of components of $(S_J\cup \rho_ZS_J) \cap H$.   Thus $\MM_J$ is not a limit leaf of the foliation $\LL'$.

It follows that $\MM_J$ is properly embedded in $\Omega:=\BB(0,1)\times (c,d)$.  For if not, $\LL'$ would have a limit leaf
in $\Omega$.  But the only leaf in $\Omega$ is $\MM_J$ itself, which, as we have just seen, is not a limit leaf.  This proves
properness.

Next we show that $\partial z/\partial \theta>0$ on $\MM_J\setminus J$. 
      (The partial derivative makes sense
because $\MM_J\setminus J$ is locally a graph.)
Since the two components of $\MM_J\setminus J$ are related by the $\rho_Z$ symmetry,
it suffices to show that $\partial z/\partial \theta>0$ on the component $M=\MM_J^+$.
Let $\nu$ be the unit normal vectorfield on $M$ given by
\[
  \nu:= \frac{(\nabla z, -1)}{W} 
\]
where $\nabla z= ( \tfrac{\partial z}{\partial x}, \tfrac{\partial z}{\partial y} )$ and $W= \left| (\nabla z, -1) \right|$.
Note that the Killing field
\[
  \frac{\partial}{\partial \theta} := (\tfrac{\partial x}{\partial\theta}, \tfrac{\partial y}{\partial \theta}, 0) = (-y,x,0)
\]
restricted to $M$ is the initial velocity vectorfield of the one-parameter 
family $\theta\mapsto R_\theta M$ of minimal surfaces.  
Thus the function
\begin{align*}
     &u: M\to \RR \\
     &u = \nu\cdot \frac{\partial}{\partial \theta}
           = W^{-1} (\tfrac{\partial z}{\partial x}\tfrac{\partial x}{\partial\theta}
           +                \tfrac{\partial z}{\partial y}\tfrac{\partial y}{\partial\theta})
           = W^{-1} \tfrac{\partial z}{\partial \theta}
\end{align*}
satisfies the Jacobi field equation
\begin{equation*}
   \Delta u + |A|^2u=0,
\end{equation*}
 where $A$ is the second fundamental form of $M$.  (See~\cite{choe}*{Lemma~1}.)

Now $\partial z/\partial \theta\ge 0$ on the disks $D_n$ 
(by~Theorem~\ref{existence2}), so $\partial z/\partial \theta\ge 0$ on $M$.  Thus $u$ is nonnegative, so 
from the Jacobi field equation, we see that $u$ is superharmonic.  
Hence by the maximum principle, $u$ is either everywhere $0$ or everywhere strictly positive.
The same is true of $\partial z/\partial \theta = W u$.
Since $\partial z/\partial \theta>0$ along the outer boundary curve $S_J$ of $M$ (by choice of $S_J$), it must
be positive everywhere on $M$.

 It remains to prove  Statement~\ref{Statement-3} of the Theorem. By Claim~\ref{FirstClaim}, the
 punctured disks  ${\mathcal D}_p$ with $p\in K$ are  leaves of $\LL'$. Each such leaf can be
 extended smoothly to include the puncture. This extension of $\LL'$ is smooth
 near any point in the interior of $K$.
 
 However, if $p\in \partial K$, then $p$ is simultaneously the limit of points in $\DD_p$
  and the limit of points in $Z\setminus K$.  At the former points, the tangent planes to $\LL'$ are
  horizontal, whereas at the latter points the tangent planes are vertical.  Thus the lamination $\LL'$
  cannot extend smoothly to any points in $\partial K$.
 \end{proof}

\section{Existence and uniqueness of embedded minimal disks with rotationally symmetric boundaries}
\label{ExistenceUniqueness}
\begin{theorem} \label{UniqueEmbeddedDisk}
Let $W$ be a nonempty, bounded, convex open subset of $\RR^3$ that is
rotationally symmetric about $Z$ and let $I= W\cap Z$. Let $R_\theta$
be rotation about $Z$ through angle $\theta$. Suppose $\Gamma$ is a
piecewise smooth simple closed curve such that $\Gamma$ contains $I$,
such that $\Gamma\setminus I$ is contained in $\partial W$, and such that
the curves
\[
    R_\theta (\Gamma\setminus \overline{I}) ,  \qquad 0\le \theta < 2\pi,
\]
 foliate $(\partial W)\setminus Z$.  Then  $\Gamma$ bounds a unique 
 embedded minimal disk $D$.  Furthermore, the disks
\[
    R_\theta D, \qquad 0\le \theta<2\pi
\]
foliate $W\setminus I$.  
\end{theorem}

\begin{proof} 
Let 
\[
   \gamma = \{ (r,z): \text{ $r\ge 0$ and $(r,0,z)\in \partial W$} \}.
\]
Let $\FF$ be the set of piecewise smooth functions $f:\gamma\to \RR$
for which the curve 
\begin{equation}\label{GammaF}
  \Gamma_f = \{ (r\cos f(r,z), r \sin f(r,z), z): (r,z)\in \gamma\} \cup I.
\end{equation}
bounds an embedded minimal disk $D_f$ whose rotated images 
   $R_\theta D_f$ foliate $W\setminus Z$.

 Note that if $f$ is constant, 
then $\Gamma _f$ is a planar curve that bounds a planar disk
whose rotated images 
foliate $W\setminus I$.   Thus $\FF$ contains all the constant functions.

\begin{claim*}
Suppose that $f\in \FF$ and that
$g:\gamma\to \RR$ is a piecewise smooth function such that
\[
    \sup \| f - g\| < \pi/4.
\]
Then $g$ must also belong to $\FF$.
\end{claim*}

Since $\FF$ is nonempty, the claim implies that
$\FF$ contains every piecewise smooth function. Thus  once we have proved the claim,
we will have proved the existence part of the theorem, because every
$\Gamma$ satisfying the hypotheses of the theorem can be written in
the form \eqref{GammaF}.

\begin{proof}[Proof of Claim]
Let
\[
   \Omega = \bigcup_{-\pi/4\le \theta \le \pi/4 }  R_\theta D_f.
\]
Note that $\Omega$ is mean convex and simply connected and that
$\Gamma_g$ is contained in  $\partial \Omega$. Hence, $\Gamma _g$
bounds a least-area disk $D_g$ in $\Omega$.

Suppose   that the $R_\theta D_g$ do not
foliate $W\setminus Z$. Then $D_g$ and $R_\theta D_g$ intersect each other
for some $\theta \in (0,\pi)$. 
Any such $\theta$ must in fact be less than $\pi/2$ because 
\begin{align*}
 D_g&\subset \Omega,   
 \\
 R_\theta D_g &\subset R_\theta \Omega, 
\end{align*}
and because $\Omega$ and $R_\theta \Omega$ are disjoint
for $\pi/2\le \theta \le \pi$.
Thus if $\alpha$ is the supremum of $\theta\in (0,\pi)$ for which $D_g$ and $R_\theta D_g$
intersect each other, then $0  < \alpha \le \pi/2$. 
The boundary curves $\Gamma_g$ and $R_\alpha\Gamma_g$ intersect only
on $\overline{I}$, so $D_g$ and $R_\alpha D_g$ must be tangent at
some point point in $D_g\cup\overline{I}$.
 At each point of $\overline{I}$,  the disks $D_g$
and $R_\alpha D_g$ make an angle of $\alpha\neq 0$ with each other.
 So the point of tangency must lie in $D_g$. This contradicts the
maximum principle.  Hence the rotated images of
$D_g$ foliate $W\setminus Z$. In particular, $D_g$ is embedded.
This completes the proof of the claim.
\end{proof}

We have proved the existence of a disk $D$ with boundary $\Gamma\cup \overline{I}$
such that the rotated images of $D$ foliate $W\setminus I$.        
It remains only to prove uniqueness.  (In this paper, we never actually use the uniqueness.)
Let
$\Sigma$ be any  embedded minimal disk with boundary $\Gamma$. Since
$\Sigma$ is embedded, it has no boundary branch points, so that
$\Sigma\cup I$ is a smooth manifold with boundary.

Since the disks $R_\theta D$ foliate $W\setminus I$,
 there is a unique continuous function 
\[
  \omega:  \overline{\Sigma} \to \RR
\]
such that
 \begin{enumerate}
 \item $\omega=0$ on $\Gamma$, \item $p\in R_{\omega(p)}D$ for $p\in
 \Sigma$, and 
 \item for $p \in I$, $R_{\omega(p)}D$ and $\Sigma$ have
 the same tangent halfplane.
 \end{enumerate}
(If $\Sigma$ were not simply connected, $\omega$ might only be
well-defined up to multiples of $2\pi$. That is, $\omega$ would take
values in $\RR/2\pi\ZZ$.  But since $\Sigma$ is simply connected, we
can lift $\omega$ to the universal cover $\RR$ of $\RR/2\pi\ZZ$.)

By the maximum principle and the boundary maximum principle
(applied to points in $I$),
the maximum value of $\omega$ must be attained on $\Gamma$.  Thus the
maximum value of $\omega$ is $0$.  Similarly the minimum value is
$0$. Thus $\omega$ is identically $0$, so $\Sigma=D$.
\end{proof}

\begin{theorem}\label{existence2}
Let
\[
    \phi: [a,b] \to \RR
\]
be a smooth, strictly increasing function. Let $\Gamma$ be the closed
curve consisting of
\[
    (\cos \phi(z), \sin \phi(z), z), \quad a \le z \le b,
\]
  together with the segment
   $I= Z \cap \{ a < z < b\}$
and two horizontal segments in the planes $z=a$ and $z=b$.

Then $\Gamma$ bounds a unique embedded minimal disk $D$. The rotated
images $R_\theta D$ foliate $C\setminus Z$, where $C$ is the cylinder 
  $\BB(0,1)\times (a,b)$.
The disk $D$ can be parametrized as
\[
   (r,\theta) \mapsto (r\cos\theta, r\sin\theta, f(r,\theta))
\]
for some function
\[
 f: [0,1] \times [\alpha, \beta] \to \RR,
\]
where $\alpha =\phi (a)$ and $\beta =\phi (b)$. That is, $D$ is a
multigraph over $\BB(0,1)\setminus \{0\}$.   Furthermore, 
$f(r,\theta)$ is a strictly increasing function of $\theta$ for each $r$. 
 In particular,  $\partial f/\partial \theta$ is everywhere nonnegative.
\end{theorem}

\begin{remark*} In fact, $\partial f/\partial \theta$ is everywhere strictly positive.
The proof
is almost identical to the proof that $\partial z/\partial \theta>0$ in
  Theorem~\ref{maintheorem}\eqref{Statement:other-leaves}.
\end{remark*}

\begin{proof} Observe that $\Gamma\setminus I$ lies on the boundary of the solid
cylinder $C= \BB(0,1)\times (a,b)$,  and that the rotated images of $\Gamma\setminus \overline{I}$ 
foliate $\partial C\setminus \overline{I}$.
Thus by Theorem~\ref{UniqueEmbeddedDisk}, $\Gamma$ bounds a unique
embedded minimal disk $D$, and
the rotated images
$R_\theta D$ foliate $C\setminus I$.

Note that since $D$ is simply connected, there is a continuous
function $\theta: \overline{D}\setminus I \to \RR$ such that for
$p=(x,y,z)\in \overline{D}\setminus I$,
\[
  (x,y) = \sqrt{x^2+y^2}\, (\cos\theta(p), \sin\theta(p)).
\]
Note also that we can choose $\theta$ so that for $(x,y,z)\in \Gamma
\setminus I$,
\[
  \theta(x,y,z) = \phi(z).
\]
 In particular, $\theta \equiv \alpha$ on $\Gamma\cap \{z=a\}$ and 
$\theta \equiv \beta$ on $\Gamma\cap \{z=b\}$, where $\alpha=\phi(a)$ and $\beta=\phi(b)$.

Since $D\cup I$ is a smooth manifold with boundary, the angle function $\theta$ extends
smoothly to $I$.  (If this is not clear, consider a point $(0,0,c)\in I$.  Then $(D\cup I)\cap \{z=c\}$,
i.e., $(D\cap \{z=c\})\cup \{(0,0,c)\}$ is a smooth curve.  Thus it has a well-defined
tangent half-line at the endpoint $(0,0,c)$, which implies that $\lim_{(x,y,c)\in D\to (0,0,c)}\theta(x,y,c)$
exists.  Define $\theta(0,0,c)$ to be that limit.  The smoothness of $D\cup I$ implies
that this extension of $\theta$ to $D\cup I$ is smooth.)

By the maximum principle, $\theta$ cannot attain its maximum or its
minimum at any interior point of $D$.  By the boundary maximum
principle, $\theta$ cannot attain its maximum or minimum at any point
of $I$.  Thus the maximum and minimum are attained on
$\Gamma\setminus I$, so the minimum value is $\alpha$ and the maximum
value is $\beta$.

To show that $D$ is a multigraph, it suffices to show that the map
\begin{equation}\label{map}
  p\in \overline{D} \mapsto ( r(p), \theta(p))  \in [0,1]\times [\alpha, \beta]
\end{equation}
is one-to-one.  (It is onto by elementary topology.)

Let 
\[
  S = \{ (p,q) \in \overline{D}\times \overline{D}: \text{ $r(p)=r(q)$ and $\theta(p)=\theta(q)$} \}.
\]
By compactness, the function
\[
  (p,q)\in S \mapsto z(q) - z(p)
\]
attains its maximum value $h=z(q_0)-z(p_0)$ at some $(p_0,q_0)\in S$.
To show that the map~\eqref{map} is one-to-one, it suffices to show that $h=0$.
To see that $h=0$, 
 let $r_0:=r(p_0)=r(q_0)$ and $\theta_0:=\theta(p_0)=\theta(q_0)$.
 
 If $r_0=1$ or if $\theta_0$ is $\alpha$ or $\beta$, then $p_0$ and $q_0$
 are both in $\Gamma\setminus I$ by the maximum principle.  But $p\mapsto (r(p),\theta(p))$ is one--to-one
 on $\Gamma\setminus I$ by choice of $\Gamma$, so in this case $h=0$ and we are done.
 
 Thus we may suppose that $r_0<1$ and that $\alpha< \theta_0<\beta$.
 Now the minimal disks $\overline{D}$ and $\overline{D}+(0,0,h)$ are tangent
 at the point $q_0$, but in neighborhood of that point $\overline{D}$ lies on one side
 of $\overline{D}+(0,0,h)$.  Thus by the strong maximum principle (if $r_0>0$) or by the 
 strong boundary maximum principle (if $r_0=0$), the two disks coincide, which
 implies that  $h=0$.
 
 It remains to show that $f(r,\theta)$ is a strictly increasing function of $\theta$
 for each $r$.
 Since the disks $R_\alpha D$ foliate $C\setminus Z$, 
it follows (for each fixed $r$) that
 the graphs of the curves
 \[ 
   \mathcal{C}_\alpha: \theta \mapsto f(r,\theta -\alpha) 
 \]
 foliate the strip $\RR\times (a,b)$.   Thus the function $\theta\mapsto f(r,\theta)$
 must be strictly monotonic.  
 Since it is strictly increasing for $r=1$,  it must be strictly increasing for all $r$.
 \end{proof}

\section{Curvature Estimates via  Rado's
Theorem}\label{RadoCurvatureSection} 
\begin{theorem}\label{RadoCurvEst}
 Let $D\subset \RR^3$ be a minimal disk contained in a vertical solid
cylinder $C=\BB\times \RR$ of radius $R$,  and let $\Gamma=\partial D$ be its boundary curve.
Suppose that $\Gamma\cap \{ a<z<
b\}$ consists of two components, each of which is a $C^1$ curve whose
tangent line has slope $\ge \eps>0$ in absolute value  at every
point.

Let
\[
    D_\delta = D \cap \{a+\delta< z < b-\delta \}
\]
where $\delta>0$. Then
\begin{enumerate}[\upshape (1)]
 \item\label{immersed} $D_\delta$ has no branch points.
 \item\label{slope}
  The slope of the tangent
plane at each point of  $D_\delta$
 is greater than or equal to
 \[
   \min \left\{ \eps, \frac{\delta}{2R} \right\}.
 \]
 \item\label{CurvatureBound} For $p\in D_\delta$,
 the norm of the second fundamental form of $D_\delta$ at $p$
 is bounded by
 \[
    \frac{C}{ \dist(p, \partial D_\delta)},
 \]
where $C$ is  constant depending only on $\eps$ and
$\delta/R$.
\end{enumerate}
\end{theorem}
The proof of Theorem~\ref{RadoCurvEst} uses the following theorem
of Rado:

\begin{theorem*}[Rado]
If the boundary of a minimal disk in $\RR^3$ intersects a plane in
fewer than four points, then $D$ has no branch points in the plane,
and $D$ intersects the plane transversely.
\end{theorem*}

See, for example,  
\cite{OssermanBook}*{Lemma 7.5} or \cite{HildebrandtBook}*{p.~272}
for a proof.

\begin{proof}[Proof of Theorem~\ref{RadoCurvEst}]
For $a<t<b$, the set  $\Gamma\cap \{z=t\}$ contains exactly two points. 
Thus by Rado's Theorem,  the surface $D\cap \{ a < z < b\}$
has no branch points and no horizontal tangent planes.

Now let $p\in D_\delta$ and let $P$ be a plane through $p$ whose
slope is less than $\eps$ and less than $\delta/(2R)$.  Since the
slope is  less than $ \delta/(2R)$, the plane $P$ does not intersect
$\BB\times (-\infty,a]$ or $\BB\times [b, \infty)$. Thus
\[
   P \cap \Gamma = P \cap(\Gamma\cap \{a<z<b\}).
\]
Since the slope of $P$ is less than $\eps$, the plane $P$ intersects each of
the the two components of $\Gamma\cap \{a<z<b\}$ in at most one
point.  Thus $P$ intersects $\Gamma$ in at most two points, so $P$ is
not tangent to $D$ at $p$ by Rado's Theorem. This proves Assertion~\ref{slope}.

Assertion~\ref{CurvatureBound} follows from 
Assertion~\ref{slope} because Assertion~\ref{CurvatureBound}
holds for any minimal surface $D_\delta$ whose image under the Gauss map omits a
nonempty open subset of the unit sphere, the constant $C$ depending
only on the size of that open set.
This curvature estimate was proved by Osserman~\cite{OssermanCurvatureBound}*{Theorem~1}, who
even obtained the optimal constant.   
Alternatively, the estimate without the optimal constant follows by a standard
blow-up argument from the fact (also due to Osserman~\cite{OssermanBook}*{Theorem~8.1})
that the image of the Gauss map of a complete, nonflat minimal surface in $\RR^3$
must be dense in the unit $2$-sphere.  
\end{proof}


\begin{bibdiv}

\begin{biblist}

\bib{Bers}{article}{
   author={Bers, Lipman},
   title={Isolated singularities of minimal surfaces},
   journal={Ann. of Math. (2)},
   volume={53},
   date={1951},
   pages={364--386},
   issn={0003-486X},
   review={\MR{0043335 (13,244c)}},
}

\bib{choe}{article}{
   author={Choe, Jaigyoung},
   title={Index, vision number and stability of complete minimal surfaces},
   journal={Arch. Rational Mech. Anal.},
   volume={109},
   date={1990},
   number={3},
   pages={195--212},
   issn={0003-9527},
   review={\MR{1025170 (91b:53007)}},
}
		

\bib{ColdingMinicozziExample}{article}{
   author={Colding, Tobias H.},
   author={Minicozzi, William P., II},
   title={Embedded minimal disks: proper versus nonproper---global versus
   local},
   journal={Trans. Amer. Math. Soc.},
   volume={356},
   date={2004},
   number={1},
   pages={283--289 (electronic)},
   issn={0002-9947},
   review={\MR{2020033 (2004k:53011)}},
}

\bib{ColdingMinicozziII}{article}{
   author={Colding, Tobias H.},
   author={Minicozzi, William P., II},
   title={The space of embedded minimal surfaces of fixed genus in a
   3-manifold. II. Multi-valued graphs in disks},
   journal={Ann. of Math. (2)},
   volume={160},
   date={2004},
   number={1},
   pages={69--92},
   issn={0003-486X},
   review={\MR{2119718 (2006a:53005)}},
   doi={10.4007/annals.2004.160.69},
}

\bib{ColdingMinicozziIV}{article}{
   author={Colding, Tobias H.},
   author={Minicozzi, William P., II},
   title={The space of embedded minimal surfaces of fixed genus in a
   3-manifold. IV. Locally simply connected},
   journal={Ann. of Math. (2)},
   volume={160},
   date={2004},
   number={2},
   pages={573--615},
   issn={0003-486X},
   review={\MR{2123933 (2006e:53013)}},
}
	
\bib{BrianDean}{article}{
   author={Dean, Brian},
   title={Embedded minimal disks with prescribed curvature blowup},
   journal={Proc. Amer. Math. Soc.},
   volume={134},
   date={2006},
   number={4},
   pages={1197--1204 (electronic)},
   issn={0002-9939},
   review={\MR{2196057 (2007d:53009)}},
}

\bib{Finn}{article}{
   author={Finn, Robert},
   title={Remarks relevant to minimal surfaces, and to surfaces of
   prescribed mean curvature},
   journal={J. Analyse Math.},
   volume={14},
   date={1965},
   pages={139--160},
   issn={0021-7670},
   review={\MR{0188909 (32 \#6337)}},
}

\bib{HildebrandtBook}{book}{
   author={Dierkes, Ulrich},
   author={Hildebrandt, Stefan},
   author={K{\"u}ster, Albrecht},
   author={Wohlrab, Ortwin},
   title={Minimal surfaces. I},
   series={Grundlehren der Mathematischen Wissenschaften [Fundamental
   Principles of Mathematical Sciences]},
   volume={295},
   note={Boundary value problems},
   publisher={Springer-Verlag},
   place={Berlin},
   date={1992},
   pages={xiv+508},
   isbn={3-540-53169-6},
   review={\MR{1215267 (94c:49001a)}},
}


\bib{SiddiqueKahn}{article}{
   author={Kahn, Siddique},
   title={A minimal lamination of the unit ball with singularities along a line segment},
   eprint={arXiv:0902.3641v2 [math.DG]},
   journal={Illinois J. of Math.},
   volume={},
   date={2008},
   status={to appear}
   number={},
   pages={},
   issn={},
   review={},
}

\bib{Kleene}{article}{
   author={Kleene, Stephen J.},
   title={A minimal lamination with Cantor set-like singularities},
   eprint={arXiv:0910.0199 [math.DG]},
   date={2009},
}

\bib{MeeksRegularity}{article}{
   author={Meeks, William H., III},
   title={Regularity of the singular set in the Colding-Minicozzi lamination
   theorem},
   journal={Duke Math. J.},
   volume={123},
   date={2004},
   number={2},
   pages={329--334},
   issn={0012-7094},
   review={\MR{2066941 (2005d:53014)}},
}

\bib{MeeksWeber}{article}{
   author={Meeks, William H., III},
   author={Weber, Matthias},
   title={Bending the helicoid},
   journal={Math. Ann.},
   volume={339},
   date={2007},
   number={4},
   pages={783--798},
   issn={0025-5831},
   review={\MR{2341900 (2008k:53020)}},
}

\bib{OssermanCurvatureBound}{article}{
   author={Osserman, Robert},
   title={On the Gauss curvature of minimal surfaces},
   journal={Trans. Amer. Math. Soc.},
   volume={96},
   date={1960},
   pages={115--128},
   issn={0002-9947},
   review={\MR{0121723 (22 \#12457)}},
}

\bib{OssermanBook}{book}{
   author={Osserman, Robert},
   title={A survey of minimal surfaces},
   edition={2},
   publisher={Dover Publications Inc.},
   place={New York},
   date={1986},
   pages={vi+207},
   isbn={0-486-64998-9},
   review={\MR{852409 (87j:53012)}},
}	


\bib{WhiteC1}{article}{
   author={White, Brian},
   title={Curvatures of embedded minimal disks
   blow up on subsets of $C^1$ curves},
   eprint={},
   journal={},
   volume={},
   date={2011},
   number={},
   pages={},
   issn={},
   review={},
   doi={},
}
	


\end{biblist}

\end{bibdiv}
\end{document}